\documentclass[10pt, letterpaper, journal, twocolumn, final]{IEEEtran}
\usepackage{amsmath}
\usepackage[hidelinks]{hyperref}
\usepackage{amssymb}
\usepackage{booktabs}
\usepackage{cite}
\usepackage{graphicx}
\usepackage{setspace}
\usepackage{epstopdf}

\allowdisplaybreaks

\renewcommand{\leq}{\leqslant}

\newcommand{\maj}{^\mathsf{M}}

\title{Convex Fused Lasso Denoising with Non-Convex Regularization and its use for Pulse Detection}
\author{
	{Ankit Parekh and Ivan W. Selesnick}
	\thanks{Copyright (c) IEEE. Personal use of this material is permitted. However, permission to use this material for any other purposes must be obtained from the IEEE by sending a request to pubs-permissions@ieee.org.}
	\thanks{A. Parekh (ankit.parekh@nyu.edu) is with the Department of Mathematics, Tandon School of Engineering, New York University, and I. Selesnick is with the Department of Electrical and Computer Engineering, Tandon School of Engineering, New York University, NY. }\thanks{This work was supported by ONR under grant N00014-15-1-2314 and by NSF under grant CCF-1525398.}\thanks{Supplementary MATLAB code is available at \url{http://goo.gl/xAi85N}}}

\begin{document}
	\maketitle
	
	%-----------------------------------------------------Abstract-----------------------------------
	\begin{abstract}
		We propose a convex formulation of the fused lasso signal approximation problem consisting of non-convex penalty functions. The fused lasso signal model aims to estimate a sparse piecewise constant signal from a noisy observation. Originally, the $\ell_1$ norm was used as a sparsity-inducing convex penalty function for the fused lasso signal approximation problem. However, the $\ell_1$ norm underestimates signal values. Non-convex sparsity-inducing penalty functions better estimate signal values. In this paper, we show how to ensure the convexity of the fused lasso signal approximation problem with non-convex penalty functions. We further derive a computationally efficient algorithm using the majorization-minimization technique. We apply the proposed fused lasso method for the detection of pulses.
	\end{abstract}
	
	
	%-----------------------------------------------------Introduction-------------------------------
	\section{Introduction}
	
	We consider the problem of estimating a sparse piecewise constant signal $x$ from its noisy observation $y$, i.e., 
	\begin{align}
		\label{eq::model}
		y = x + w, \qquad y,x,w \in \mathbb{R}^N,
	\end{align} where $w$ represents zero-mean additive white Gaussian noise. Estimation of sparse piecewise continuous signals arise in transient removal \cite{Selesnick2014_TARA}, genomic hybridization \cite{Rapaport2008,Tibshirani2008}, signal and image denoising \cite{Friedman2007,Bayram2014}, prostate cancer analysis \cite{Tibshirani2005}, sparse trend filtering \cite{Tibshirani2014,Ramdas2014,Kim2009,Wang2014} and biophysics \cite{Little2011}. In order to estimate sparse piecewise constant signals, it has been proposed \cite{Tibshirani2005} to solve the following sparse-regularized optimization problem 
	\begin{align}
		\label{eq:: L1 fused lasso}
		\arg\min_{x} \biggl\lbrace \dfrac{1}{2}\|y-x\|_2^2 + \lambda_0\|x\|_1 + \lambda_1 \|Dx\|_1 \biggr\rbrace,
	\end{align}where $\lambda_0 > 0$ and $\lambda_1 > 0$ are the regularization parameters and the matrix $D$ is defined as
	\begin{align}
		\label{eq:matrix D}
		D = \left[ \begin{array}{cccc}
			-1 & 1 & & \\
			& \ddots & \ddots & \\
			& & -1 & 1
		\end{array} \right], \quad D \in \mathbb{R}^{(N-1) \times N}.
	\end{align} The optimization problem \eqref{eq:: L1 fused lasso} is well-known as the $\ell_1$ fused lasso signal approximation (FLSA) problem \cite{Tibshirani2005}. The FLSA problem has been explored to aid the diagnosis of Alzheimer's disease \cite{Xin2014} and background subtraction problems \cite{Xin2015}, however with a different data-fidelity term. Note that when $\lambda_0 = 0$, problem \eqref{eq:: L1 fused lasso} reduces to the total variation denoising (TVD) problem \cite{Rudin1992}. 
	
	It is known that the $\ell_1$ norm, when used as a sparsity-inducing regularizer, underestimates the signal values. The $\ell_1$ norm is generally not the tightest convex envelope for sparsity \cite{Jojic2011}. In order to better estimate signal values, non-convex penalty functions are often favored over the $\ell_1$ norm \cite{Bruckstein2009,Candex_2008_JFAP,Chartrand2015,Portilla2007,Nikolova2011,Repetti2014}. However, the use of non-convex penalty functions generally leads to non-convex optimization, which suffer from several issues (spurious local minima, initialization, convergence, etc.). 
	
	In this paper, we propose to estimate sparse piecewise constant signals via the following convex non-convex (CNC) FLSA problem 
	\begin{multline}
		\label{eq::cost function}
		\arg\min_{x} \biggl\lbrace F(x) = \dfrac{1}{2}\|y-x\|_2^2 + \lambda_0\sum_{n=1}^{N} \phi(x_n;a_0) \\
		+ \lambda_1\sum_{n=1}^{N-1} \phi([Dx]_n;a_1)\biggr\rbrace, 
	\end{multline} where $\phi\colon \mathbb{R} \to \mathbb{R}$ with $a\geqslant 0$ is a non-convex sparsity-inducing regularizer. Specifically, we propose that the regularization terms be chosen so that the objective function $F$ in \eqref{eq::cost function} is convex. As a result, the CNC FLSA approach avoids the drawbacks of non-convex optimization. The non-convex penalty function $\phi$ aims to induce sparsity more strongly than the $\ell_1$ norm and thus better estimate the signal values. The parameter $a$ controls the degree of non-convexity of $\phi$; a higher value of $a$ indicates a higher degree of non-convexity for $\phi$. As a main result, we state and prove a condition that $a_0$ and $a_1$ must satisfy to ensure the objective function $F$ in \eqref{eq::cost function} is strictly convex. As a consequence of the convexity condition, well-known convex optimization techniques can be used to reliably obtain the global minimum of the objective function $F$. As a second main result, we provide an efficient fast converging algorithm for the proposed CNC FLSA problem \eqref{eq::cost function}. The algorithm is derived using the majorization-minimization (MM) procedure \cite{Figueiredo2007}. 
	
	The idea of formulating a convex problem with non-convex regularization was described by Blake and Zisserman \cite{Blake1987}, and Nikolova \cite{Nikolova1998,Nikolova1999}. The idea is to balance the positive second-derivatives of the data-fidelity term with the negative second-derivatives of the non-convex penalty function. This approach has been successfully applied to various signal processing applications (e.g., \cite{Selesnick2015,Parekh2015,Lanza2015,Selesnick2015_2} and the references therein). Using this technique, a modified formulation of the $\ell_1$ FLSA problem \eqref{eq:: L1 fused lasso} was proposed with an aim to induce sparsity more strongly than the $\ell_1$ norm \cite{Bayram2014}. Further, a two-step procedure was used to obtain the solution to the modified fused lasso problem. However, the modified fused lasso (MDFL) problem \cite{Bayram2014} considers only the first of the two regularization terms as non-convex. 
	
	This paper is organized as follows. In Section~\ref{section::Preliminaries} we describe the class of non-convex penalty functions. In Section~\ref{section::Convexity condition} we  provide the convexity condition for the objective function $F$ in \eqref{eq::cost function}. We derive an algorithm to solve the CNC FLSA problem \eqref{eq::cost function} based on the MM procedure in Section~\ref{section::Algorithm}. In Section~\ref{section::Examples} we apply the proposed CNC FLSA approach to the problem of detecting ECG pulses in strong additive white Gaussian noise (AWGN).

	\section{Preliminaries}
	\label{section::Preliminaries}
	\subsection{Notation}
	We denote vectors and matrices by lower and upper case letters respectively. The $N$-point signal $y$ is represented by the vector
	\begin{align}
		y = [y_0, \hdots, y_{N-1}]^{T}, \qquad y \in \mathbb{R}^N,
	\end{align} where $[\cdot]^T$ represents the transpose. The $\ell_1$ and $\ell_2$ norms of the vector $y$ are defined as
	\begin{align}
		\|y\|_1 = \sum_n|y(n)|, \quad \|y\|_2 = \left(\sum_n|y(n)|^2\right)^{1/2}.
	\end{align} The soft-threshold function \cite{Donoho1995} for $\lambda>0, \lambda \in \mathbb{R}$ is defined as
	\begin{align}
		\label{eq::softThresholdFunction}
		\mbox{soft}(x;\lambda) = \begin{cases}
			x +\lambda, & \hfil x < -\lambda\\
			\hfil 0, & \hfil -\lambda \leqslant x \leqslant \lambda\\
			x-\lambda, & \hfil x > \lambda.
		\end{cases}
	\end{align}For $x \in \mathbb{R}^N$, the notation $\mbox{soft}(x;\lambda)$ implies that the soft-threshold function is applied element-wise to $x$ with a threshold of $\lambda$.
	\newtheorem{Definition}{\bf Definition}
	\begin{Definition}
		The total variation denoising (TVD) problem \cite{Rudin1992} is defined as
		\begin{align}
			\label{eq::TV denoising}
			\mbox{tvd}(y;\lambda) = \arg\min_{x} \left\lbrace\dfrac{1}{2}\|y-x\|_2^2 + \lambda \|Dx\|_1 \right\rbrace,
		\end{align}where $\lambda > 0$ is the regularization parameter. 
	\end{Definition}

	We note the following lemma, which provides an efficient two-step solution to the $\ell_1$ FLSA problem \eqref{eq:: L1 fused lasso}. 
	\newtheorem{lemma}{\bf Lemma}
	\begin{lemma}{\cite[Lemma A.1]{Friedman2007}}
		The solution $x^*$ to the $\ell_1$ FLSA problem \eqref{eq:: L1 fused lasso} is given by 
		\begin{align}
			\label{eq::Solution to L1 fused lasso}
			x^* = \mbox{soft}\bigl(\mbox{tvd}(y,\lambda_1),\lambda_0\bigr).
		\end{align}
	\end{lemma}

	\subsection{Non-convex penalty functions}
	\label{subsec::Penalty}
	We propose to use parameterized non-convex penalty functions, with a view to induce sparsity more strongly than the $\ell_1$ norm. We assume such non-convex penalty functions have the following properties. 
	
	\newtheorem{assumption}{\bf Assumption}
	\begin{assumption}
		\label{theorem::assumption 1}
		The penalty function $\phi\colon\mathbb{R} \to \mathbb{R}$ satisfies the following
		\begin{enumerate}
			\item $\phi$ is continuous on $\mathbb{R}$, twice differentiable on $\mathbb{R}\!\setminus\! \lbrace 0\rbrace$ and symmetric, i.e., $\phi(-x; a) = \phi(x; a)$
			\item $\phi'(x) > 0, \forall x > 0$	
			\item $\phi''(x) \leq 0, \forall x > 0$
			\item $\phi'(0^{+}) = 1$
			\item $\inf\limits_{x\neq0}\phi''(x;a) = \phi''(0^+;a) = -a$
		\end{enumerate}
	\end{assumption}
	
	An example of a penalty function, which satisfies Assumption \ref{theorem::assumption 1}, is the logarithmic penalty function \cite{Candex_2008_JFAP} defined as
	\begin{align}
		\label{eq::log penalty}
		\phi(x;a) = \begin{cases}
			\dfrac{1}{a}\log(1 + a|x|), & a > 0\\
			|x|, & a = 0.
		\end{cases}
	\end{align}Note that when $a=0$, this penalty function reduces to the $\ell_1$ norm. Other examples of non-convex penalty functions satisfying Assumption \ref{theorem::assumption 1} include the arctangent and the rational penalty functions \cite{Geman1992,Selesnick_2014_MSC}. Note that the $\ell_p$ norm does not satisfy Assumption \ref{theorem::assumption 1}.

	We note the following lemma, which we will use to obtain a convexity condition for optimization problem \eqref{eq::cost function}.
	\begin{lemma}{\cite{Parekh2015}}
		\label{theorem::Lemma 2}
		Let $\phi\colon\mathbb{R}\to\mathbb{R}$ satisfy Assumption \ref{theorem::assumption 1}. The function $s\colon\mathbb{R}\to\mathbb{R}$ defined as
		\begin{align}
			\label{eq::function s}
			s(x;a)
			= \phi(x;a) - |x|,
		\end{align}is twice continuously differentiable and concave with
		\begin{align}
			-a \leq s''(x;a) \leq 0.
		\end{align}
	\end{lemma}

	\section{Convexity condition}
	\label{section::Convexity condition}
	In this section, we seek to find a condition on the parameters $a_0$ and $a_1$ to ensure that the objective function $F$ in \eqref{eq::cost function} is strictly convex. The following theorem provides the required condition on $a_0$ and $a_1$. 
	\newtheorem{theorem}{\bf Theorem}
	\begin{theorem}
		\label{theorem::Theorem 1}
		Let $\phi\colon\mathbb{R}\to\mathbb{R}$ be a non-convex penalty function satisfying Assumption 1. The function $F\colon\mathbb{R}^N\to\mathbb{R}$ defined as
		\begin{align}
			\label{eq::definition of F}
			F(x) = \dfrac{1}{2}\|y-x\|_2^2 + \lambda_0\sum_{n=1}^{N} \phi(x_n;a_0) + \lambda_1\sum_{n=1}^{N-1} \phi([Dx]_n;a_1),
		\end{align} is strictly convex if
		\begin{align}
			\label{eq::condition on a1, a2}
			0 \leqslant a_0\lambda_0 + 4a_1\lambda_1 \leqslant 1.
		\end{align}
	\end{theorem}
	
	\begin{IEEEproof}
		Consider the function $G:\mathbb{R}^N \to \mathbb{R}$ defined as
		\begin{align}
			\label{surrogate cost function}
			G(x) = \dfrac{1}{2}\|y-x\|_2^2 + \lambda_0 \sum_{n=1}^{N} s(x_n;a_0) + \lambda_1 \sum_{n=1}^{N-1} s([Dx]_n;a_1),
		\end{align}where $s(x;a) = \phi(x;a) - |x|$. From Lemma \ref{theorem::Lemma 2}, the function $G$ is twice continuously differentiable and its Hessian can be written as
		\begin{align}
			\label{eq::Hessian of G}
			\nabla^2 G = I + \lambda_0\Gamma(x;a_0) + \lambda_1 D^{T} \Gamma(Dx;a_1) D,
		\end{align} where
		\begin{align}
			\Gamma(x;a) = \left[ \begin{array}{ccc}
				s''_1(x_1;a) & & \\
				& \ddots & \\
				& & s''_{n}(x_n;a)
			\end{array}\right].
		\end{align} For the strict convexity of $G$, we need to ensure that $\nabla^2 G$ is positive definite. To this end, from the assumptions on $\phi$, it follows that 
		\begin{align}
			\label{eq::Inequality for a0}
			\Gamma (x;a_0) &\succcurlyeq -a_0 I, \quad x \in \mathbb{R}^N.\\
			\intertext{Moreover, we can write}
			D^T\Gamma (Dx;a_1) D &\succcurlyeq -a_1 D^TD \\
			\label{eq::inequality using eigenvalues}
			&\succ -4a_1 I.
		\end{align} The inequality \eqref{eq::inequality using eigenvalues} is obtained using the eigenvalues\footnote{The eigenvalues of $D^TD$ are given by $\lbrace 2 - 2\cos(k\pi/N)\rbrace$ for $k = 0,\hdots,N-1$ \cite{Strang1999}.} of the matrix $D^TD$. Using \eqref{eq::Hessian of G}, \eqref{eq::Inequality for a0} and \eqref{eq::inequality using eigenvalues}, $\nabla^2 G \succ 0$ if
		\begin{align}
			(1 - a_0\lambda_0 -4a_1\lambda_1) I \succcurlyeq 0,
			\intertext{or equivalently if,}
			\label{eq::inequality on a1 a2}
			1 - a_0\lambda_0 -4a_1\lambda_1 \geqslant 0.
		\end{align} From \eqref{eq::function s}, \eqref{eq::definition of F} and \eqref{surrogate cost function} it is straighforward that
		\begin{align}
			F(x) = G(x) + \lambda_0\|x\|_1 + \lambda_1 \|Dx\|_1.
		\end{align} Hence, $F$ in \eqref{eq::definition of F} is strictly convex as long as the inequality \eqref{eq::inequality on a1 a2} holds true (the function $F$ is a sum of a strictly convex function $G$, the convex $\ell_1$ norm, and the convex TV penalty). 
	\end{IEEEproof}
	
	\begin{figure}
		\centering
		\includegraphics[]{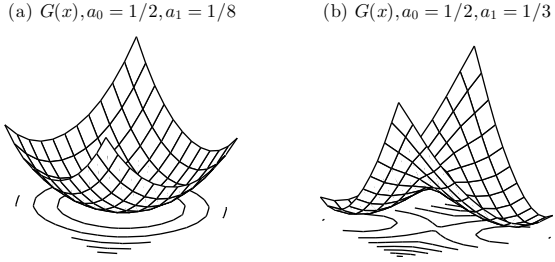}
		\caption{Surface plots illustrating the convexity condition. (a) The function $G(x)$ \eqref{surrogate cost function} is convex for $a_0 = 1/2$ and $a_1 = 1/8$. (b) The function $G(x)$ is not convex for $a_0 = 1/2$ and $a_1 = 1/3$ (these values violate Theorem 1).}
		\label{surface plots}
	\end{figure}
	
	The following example illustrates the convexity condition \eqref{eq::inequality on a1 a2} for $N = 2$. Let $\lambda_0 = \lambda_1 = 1$ and $y=0$. As per Theorem 1, the function $G$ (by extension the function $F$) is strictly convex if $a_0 + 4a_1 \leqslant 1$. Figure \ref{surface plots}(a) shows the function $G$ with the values $a_0 = 1/2$ and $a_1 = 1/8$. These values satisfy Theorem 1 and as a result the function $G$ is strictly convex. On the other hand, Fig.~\ref{surface plots}(b) shows the function $G$ when $a_0 = 1/2$ and $a_1 = 1/3$. These values of $a_0$ and $a_1$ violate the Theorem~1; consequently the function $G$ is non-convex as seen in Fig.~\ref{surface plots}(b).

	\begin{figure}
		\centering
		\includegraphics[]{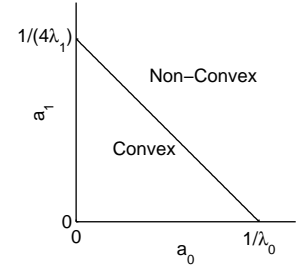}
		\caption{Region of convexity for the function $F$ in \eqref{eq::definition of F}. The function $F$ is strictly convex for any values of $a_0$ and $a_1$ inside the triangular region.}
		\label{fig::Convexity region}
	\end{figure}
	The convexity condition given by Theorem 1 in \eqref{eq::condition on a1, a2} implies that the values of $a_0$ and $a_1$ must lie on or below the line given by $a_0\lambda_1 + 4a_1\lambda_1 = 1$. Figure \ref{fig::Convexity region} displays the values of $a_0$ and $a_1$ for which the function $F$ is strictly convex. In order to maximally induce sparsity, we choose the values of $a_0$ and $a_1$ on the line. Specifically, we propose to select a value of $a_0 \in (0, 1/\lambda)$ and set the value of $a_1$ as 
	\begin{align}
		\label{eq::a_1 value}
		a_1 = \dfrac{1-a_0\lambda_0}{4\lambda_1}.
	\end{align}

	\section{Optimization Algorithm}
	\label{section::Algorithm}
	Due to Theorem \ref{theorem::Theorem 1}, we can reliably obtain via convex optimization the global minimum of \eqref{eq::cost function} as long as the parameters $a_0$ and $a_1$ are chosen to satisfy \eqref{eq::condition on a1, a2}. We derive an algorithm for the proposed CNC fused lasso method using the majorization-minimization (MM) procedure \cite{Figueiredo2007}, such that 
	\begin{align}
		\label{eq::first mm update}
		x^{k+1} &= \arg\min_{x} F\maj (x,x^{k}),
	\end{align} where $F\maj$ denotes a majorizer of the function $F$ in \eqref{eq::cost function}, and where $k$ is the iteration index. The MM procedure guarantees that each iteration monotonically decreases the value of the objective function $F$ in \eqref{eq::cost function}. We use the absolute value function and a linear function to majorize the non-convex penalty function. With this particular choice of majorizer, each MM update iteration involves solving the $\ell_1$ FLSA problem \eqref{eq:: L1 fused lasso}.
	
	\begin{figure}
		\centering
		\includegraphics[]{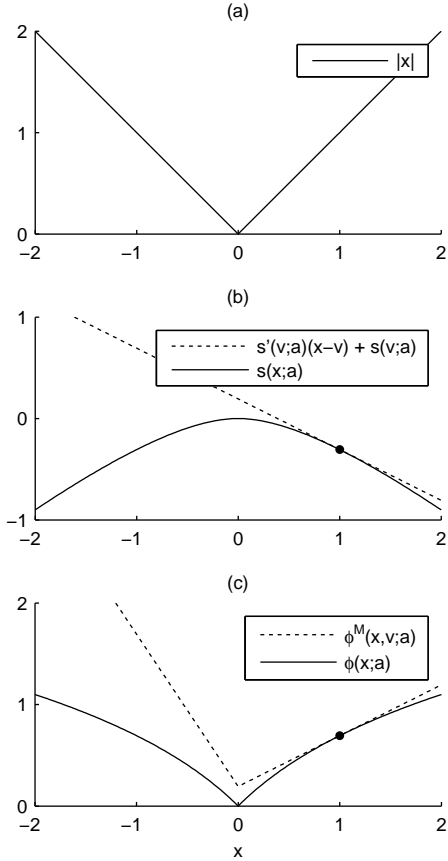}
		\caption{(a) The absolute value function $|x|$. (b) The twice continuously differentiable function $s(x;a)$ and the tangent line at $x=1$. (c) The non-convex penalty function $\phi$ and its majorizer $\phi\maj(x,v;a)$ given in \eqref{eq:: Majorizer of phi}.}
		\label{fig::Majorizer}
	\end{figure}
	
	To derive a majorizer of the function $\phi$, note that $\phi(x;a) = s(x;a) + |x|$. As a result, it suffices to majorize the function $s$ with a linear term in order to obtain a majorizer of the function $\phi$. Observe that since $s$ is a concave function, the tangent line to $s$ at a point $v$ always lies above the function $s$. Using the tangent line to the function $s$, a majorizer of the function $\phi$ is given by $\phi\maj \colon\mathbb{R}\times \mathbb{R} \to \mathbb{R}$, defined as
	\begin{align}
		\label{eq:: Majorizer of phi}
		\phi\maj (x,v;a) = |x| + s'(v;a)(x-v) + s(v;a),
	\end{align}for $ x,v \in \mathbb{R}$. It follows straightforwardly that
	\begin{align}
		\label{eq::Conditions of majorizer}
		\phi\maj (x,v;a) &\geqslant \phi(x;a), \quad \forall x,v \in \mathbb{R}, \\
		\label{eq::Conditions of majorizer 2}
		\phi\maj (v,v;a) & = \phi(v;a), \quad \forall v \in \mathbb{R}.
	\end{align}
	
	Figure~\ref{fig::Majorizer}(a) shows the absolute value function $|x|$. The twice continuously differentiable function $s(x;a)$ is shown in Fig.~\ref{fig::Majorizer}(b), along with the tangent line to $s(x;a)$ at $x=1$. Figure~\ref{fig::Majorizer}(c) shows the non-convex penalty function $\phi$ and its majorizer $\phi\maj(x,v;s)$ given by \eqref{eq:: Majorizer of phi}. The majorizer is the sum of the absolute value function in Fig.~\ref{fig::Majorizer}(a) and the tangent line to $s(x;a)$ in Fig.~\ref{fig::Majorizer}(b).
	
	Using \eqref{eq::Conditions of majorizer} and \eqref{eq::Conditions of majorizer 2}, we note that
	\begin{align}
		\label{eq::Majorize condition on Dx and x}
		\sum_{n} \phi\maj \bigl(x_n,v_n;a \bigr) &\geqslant \sum_n \phi\bigl(x_n;a\bigr), \\
		\sum_{n} \phi\maj \bigl([Dx]_n,[Dv]_n;a \bigr) &\geqslant \sum_n \phi\bigl([Dx]_n;a\bigr), 
	\end{align}with equality if $x =v$. Further, note that
	\begin{align}
		\label{eq::Simplification of majorizers a}
		\sum_n \phi\maj \bigl(x_n,v_n;a \bigr) &= \|x\|_1 + s'(v;a)^T (x-v) + C_1,
	\end{align}where $s'(v;a)$ is the vector defined as $[s'(v;a)]_n = s'(v_n;a)$, i.e., the derivative of the function $s$ is applied element-wise to the vector $v$. Further, note that $C_1$ is a constant that does not depend on $x$. Similarly, we write 
	\begin{align}
		\label{eq::Simplification of majorizers b}
		\sum_{n} &\phi\maj \bigl([Dx]_n,[Dv]_n;a \bigr) \nonumber \\
		&\qquad = \|Dx\|_1 + s'\bigl(Dv;a\bigr)^TD(x-v) + C_2,
	\end{align} where $C_2$ is a constant that does not depend on $x$. Therefore, using \eqref{eq::Simplification of majorizers a} and \eqref{eq::Simplification of majorizers b}, a majorizer of the objective function $F$ in \eqref{eq::cost function} is given by $F\maj \colon\mathbb{R}^N \times \mathbb{R}^N \to \mathbb{R}$, defined as
	\begin{align}
		\label{eq::Majorizer of F}
		F\maj (x,v) &= \dfrac{1}{2}\|y-x\|_2^2 + \lambda_0\|x\|_1 + \lambda_0 s'(v;a)^T(x-v) \nonumber \\
		&\quad+ \lambda_1\|Dx\|_1 + \lambda_1  s'\bigl( Dv;a\bigr)^TD(x-v) + C,
	\end{align} where $C$ is a constant that does not depend on $x$. Completing the square, we write \eqref{eq::Majorizer of F} as
	\begin{align}
		\label{eq::Simplified majorizer of F}
		F\maj (x,v) &= \dfrac{1}{2} \| \tilde{y}(v) - x\|_2^2 + \lambda_0\|x\|_1 + \lambda_1\|Dx\|_1 + C,
	\end{align}where  
	\begin{align}
		\label{eq::tilde y}
		\tilde{y}(v) = y - \lambda_0 s'(v;a) - \lambda_1 D^T s'\bigl(Dv;a \bigr).
	\end{align}Therefore, each MM iteration consists of minimizing the function $F\maj$ \eqref{eq::Simplified majorizer of F}, which is the $\ell_1$ FLSA problem \eqref{eq:: L1 fused lasso} with $\tilde{y}(v)$ as the input. Consequently, using \eqref{eq::Solution to L1 fused lasso} the MM update \eqref{eq::first mm update} can be written as
	\begin{subequations}
		\label{eq:: MM algorithm}
		\begin{align}
			\tilde{y}^{k} &= y - \lambda_0 s'(x^{k};a) - \lambda_1 D^T s'\bigl(Dx^{k};a \bigr) \\
			\label{eq:: MM update iteration}
			x^{k+1}	&= \mbox{soft}\bigl(\mbox{tvd(}\tilde{y}^{k},\lambda_1),\lambda_0\bigr).
		\end{align}
	\end{subequations} 
	
	Equation \eqref{eq:: MM algorithm} constitutes a fast converging, computationally efficient algorithm to solve the proposed CNC FLSA problem \eqref{eq::cost function}. To implement \eqref{eq:: MM update iteration}, we use the fast (finite-time) exact TV denoising algorithm based on the `taut-string' method \cite{Condat2013}, which has a worst case complexity of $\mathcal{O}(n)$. We initialize the iteration with the solution \eqref{eq::Solution to L1 fused lasso} to the $\ell_1$ FLSA problem \eqref{eq:: L1 fused lasso}. Note that the MM update \eqref{eq:: MM algorithm} does not consist of any matrix inverses.

	\section{Examples}
	\begin{figure}
		\centering
		\includegraphics[]{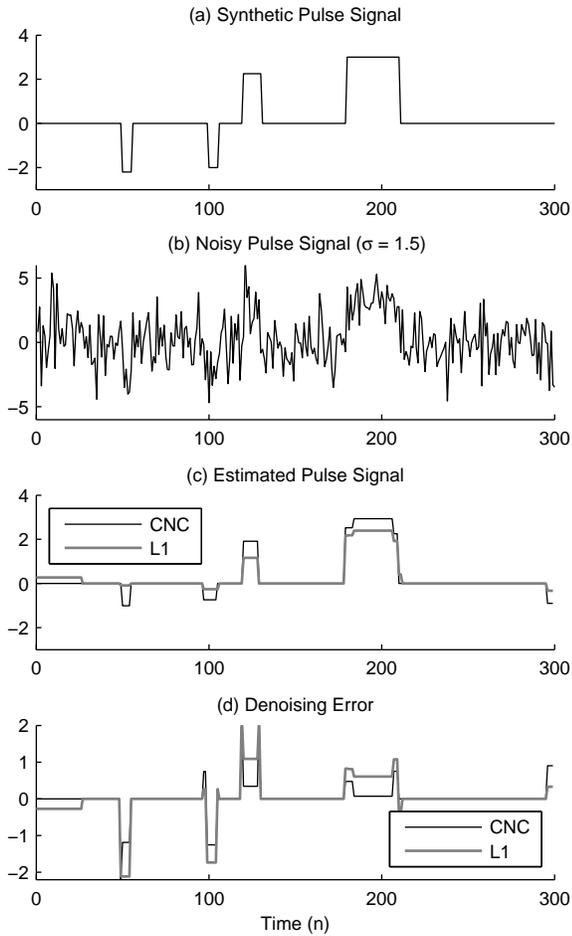}
		\caption{Denoising a synthetic signal using CNC fused lasso \eqref{eq::cost function} and $\ell_1$ fused lasso \eqref{eq:: L1 fused lasso}.}
		\label{fig::pulse signal denoising}
	\end{figure}
	\label{section::Examples}
	We consider the problem of estimating pulses of varying width in the presence of high additive white Gaussian noise (AWGN). We model the pulse signal as sparse piecewise constant and apply the proposed CNC FLSA problem \eqref{eq::cost function} to estimate the individual pulses. We set the value of $\lambda_1$ as in \cite{Selesnick2015}, i.e., $\lambda_1 = \beta\sqrt{N}\sigma$ where $\beta$ is a constant (usually $1/4$) and $\sigma$ represents the standard deviation of the additive white Gaussian noise. We manually set the value of $\lambda_0$ to obtain the lowest RMSE.

	\begin{figure}[t!]
		\centering
		\includegraphics[]{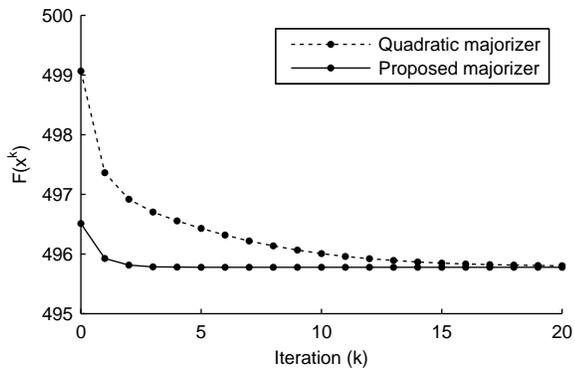}
		\caption{The value of the objective function $F$ in \eqref{eq::cost function} at each iteration of the MM algorithm using quadratic and proposed majorizers. The proposed MM algorithm \eqref{eq:: MM algorithm} converges within 5 iterations.}
		\label{fig::Cost function history}
	\end{figure}
	
	Figure~\ref{fig::pulse signal denoising}(a) illustrates the synthetic clean pulse signal and Fig.~\ref{fig::pulse signal denoising}(b) shows the noisy pulse signal. Shown in Fig.~\ref{fig::pulse signal denoising}(c) are the estimates obtained using the standard $\ell_1$ norm and the non-convex atan penalty function [19, equation (23)]. It can be seen that the proposed CNC FLSA method estimates the pulses more accurately than the $\ell_1$ FLSA method. The relative performance of the CNC fused lasso in estimating pulses is also highlighted by the denoising error shown in Fig.~\ref{fig::pulse signal denoising}(d). 
	
	The value of the objective function $F$ in \eqref{eq::cost function}, after each iteration of the MM algorithm \eqref{eq:: MM algorithm}, is shown in Fig.~\ref{fig::Cost function history}. The MM algorithm derived in \cite[Table~II]{Selesnick2014_TARA} for a more general problem can also be used to solve the CNC FLSA problem \eqref{eq::cost function}. However, the MM algorithm in \cite{Selesnick2014_TARA} utilizes a quadratic majorizer for the non-convex penalty function $\phi$. Figure~\ref{fig::Cost function history} shows that the proposed MM algorithm \eqref{eq:: MM algorithm} converges much faster than the MM algorithm in \cite{Selesnick2014_TARA}. The proposed MM algorithm \eqref{eq:: MM algorithm} converges in about 5 iterations. 
	
	\begin{figure}
		\centering
		\includegraphics[]{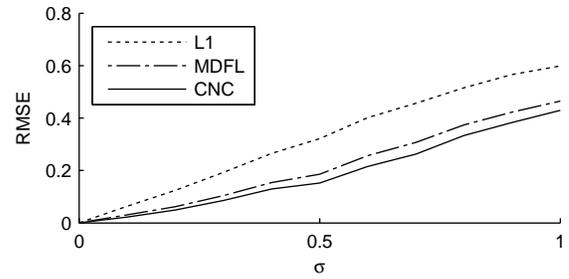}
		\caption{Average RMSE as a function of $\sigma$. The proposed CNC FLSA yields the lowest RMSE across different values of the noise variance $(\sigma^2)$.}
		\label{fig::average rMSE}
	\end{figure}
	
	\begin{figure}
		\centering
		\includegraphics[]{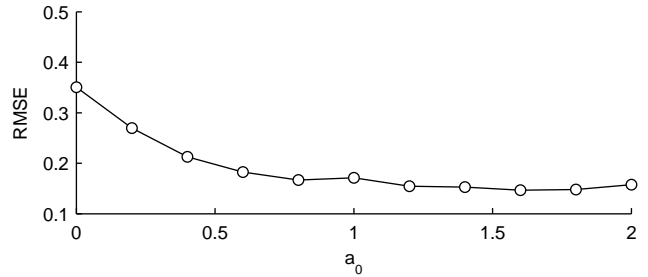}
		\caption{Average RMSE as a function of $a_0$ for the synthetic test signal example in Fig.~\ref{fig::pulse signal denoising}. }
		\label{fig::RMSE as a function of a0}
	\end{figure}

	In order to assess the relative performance of the proposed CNC fused-lasso method, we use 15 realizations of the noisy synthetic pulse signal in Fig.~\ref{fig::pulse signal denoising} and denoise them using both the original $\ell_1$ fused lasso and the proposed CNC fused lasso methods. We also compare with the modified fused lasso (MDFL) \cite{Bayram2014}, which is a special case of \eqref{eq::cost function} with $a_1 = 0$; i.e., only the first regularization term is non-convex. It can be seen in Fig.~\ref{fig::average rMSE} that the proposed CNC FLSA \eqref{eq::cost function} approach offers the lowest RMSE values across different noise levels. Further, for the test signal in Fig.~\ref{fig::pulse signal denoising}(a), the average RMSE as a function of $a_0$ is shown in Fig.~\ref{fig::RMSE as a function of a0}. Note that $a_1$ is set according to \eqref{eq::a_1 value}.

	\begin{figure}[t!]
		\centering
		\includegraphics{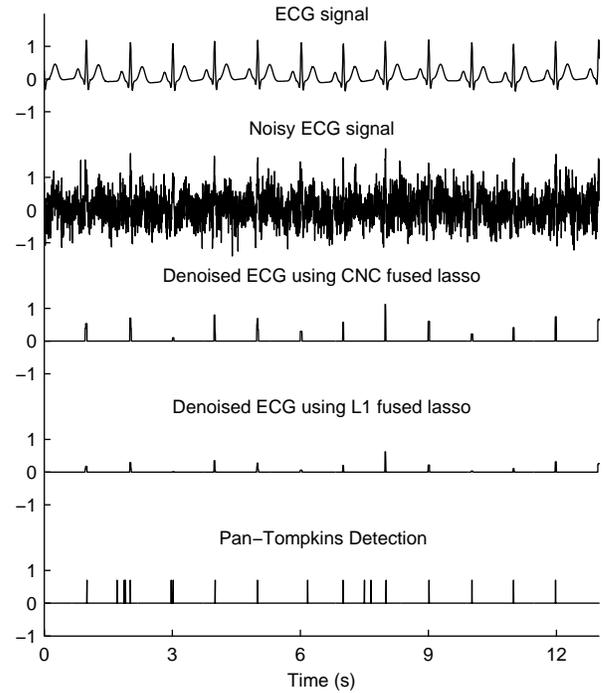}
		\caption{Denoising of ECG signal in strong AWGN. The $\ell_1$ FLSA underestimates the signal values. The Pan-Tompkins detects several false-positive R-waves. }
		\label{fig::ECG detection}
	\end{figure}

	As an another example, we consider the problem of detecting the QRS peaks in an ECG signal in AWGN with high variance ($\sigma^2$). Wearable heart-rate monitors suffer from strong noise due to abrupt motion artifacts. Several methods for detecting the QRS peaks in ECG signals were studied in \cite{Friesen1990,Pan1985}. We evaluate the detection of ECG R-waves in strong AWGN using the proposed CNC fused lasso method. Using a sampling frequency of 256 Hz, we simulate the ECG signal using the synthetic ECG waveform generator, ECGSYN \cite{McSharry2003}. The clean and noisy ECG ($\sigma = 0.4$) are shown in Fig.~\ref{fig::ECG detection}.  We set the parameters $\lambda_0 = 0.6$, $\lambda_1 = 0.9$, $a_0 = 0.9/\lambda_0$ and $a_1 = 0.1/(4\lambda_1)$. We use 20 iterations for the proposed CNC FLSA algorithm \eqref{eq:: MM algorithm}. 
	
	Figure~\ref{fig::ECG detection} illustrates the denoised ECG signal using the $\ell_1$ FLSA and the proposed CNC FLSA methods. It can be seen that the $\ell_1$ FLSA does not detect all the R-waves. Moreover, the amplitudes of the pulses detected using the proposed CNC FLSA \eqref{eq::cost function} are relatively high compared to those detected using the $\ell_1$ FLSA. The $\ell_1$ norm tends to underestimate signal values. Also shown in Fig.~\ref{fig::ECG detection} are the R-waves detected using the Pan-Tompkins real-time QRS detector \cite{Pan1985}. Note that the Pan-Tompkins detector was not designed for ECG signals with high-noise variance. As a result, the Pan-Tompkins algorithm detects several false-positive R-waves. 
	
	\section{Conclusion}
	
	The fused lasso signal approximation (FLSA) problem aims to estimate sparse piecewise constant signals. In order to improve the accuracy of the $\ell_1$ FLSA approach, we use non-convex penalty functions as sparsity-inducing regularizers. In this paper we generalize the results of \cite{Bayram2014}, which addresses the case wherein only one of the two regularization terms is non-convex. We prove that the proposed FLSA objective function is convex when the non-convex penalty parameters are suitably set. We also derive a computationally efficient algorithm using the majorization-minimization technique. The proposed CNC FLSA algorithm does not consist of any calculations involving a matrix inverse. We apply the proposed method to the problem of pulse detection under high additive white Gaussian noise. An illustration is provided for the detection of R-waves in an ECG signal. 
	
	\section{Acknowedgement}
	The authors would like to thank Ilker Bayram for valuable comments on an earlier version of this paper. 
	
	% Generated by IEEEtranS.bst, version: 1.13 (2008/09/30)

\end{document}